\documentclass[12pt]{article}
\usepackage{amsfonts}
\usepackage{amssymb}
\usepackage{enumerate}

\newtheorem{ttt}{Theorem}[section]
\newtheorem{llll}[ttt]{Lemma}
\newtheorem{ccc}[ttt]{Claim}
\newtheorem{eee}[ttt]{Example}
\newtheorem{sss}[ttt]{Statement}
\newtheorem{ddd}[ttt]{Definition}
\newtheorem{qqq}[ttt]{Question}
\newtheorem{cccc}[ttt]{Corollary}

\newcommand{\bt}{\begin{ttt}}
\newcommand{\bl}{\begin{llll}}
\newcommand{\bc}{\begin{ccc}}
\newcommand{\bex}{\begin{eee}}
\newcommand{\bs}{\begin{sss}}
\newcommand{\bd}{\begin{ddd} \upshape}
\newcommand{\bq}{\begin{qqq}}
\newcommand{\bcor}{\begin{cccc}}

\newcommand{\bp}{\noindent\textbf{Proof }}
\newcommand{\br}{\noindent\textbf{Remark }}
\newcommand{\ba}{\bigskip\noindent\textbf{Acknowledgements }}

\newcommand{\et}{\end{ttt}}
\newcommand{\el}{\end{llll}}
\newcommand{\ec}{\end{ccc}}
\newcommand{\eex}{\end{eee}}
\newcommand{\es}{\end{sss}}
\newcommand{\ed}{\end{ddd}}
\newcommand{\eq}{\end{qqq}}
\newcommand{\ecor}{\end{cccc}}

\newcommand{\ep}{\hspace{\stretch{1}}$\square$\medskip}
\newcommand{\er}{\medskip}

\newcommand{\lab}[1]{\label{#1}}

\newcommand{\NN}{\mathbb{N}}

\newcommand{\RR}{\mathbb{R}}

\newcommand{\al}{\alpha}

\newcommand{\om}{\omega}

\newcommand{\ka}{\kappa}
\newcommand{\la}{\lambda}

\newcommand{\Hd}{\mathcal{H}^d} 
\newcommand{\Ss}{Sierpi\'nski set}
\newcommand{\Sss}{Sierpi\'nski sets}  
\newcommand{\Ssc}{Sierpi\'nski set in the sense of cardinality} 
\newcommand{\Sssc}{Sierpi\'nski sets in the sense of cardinality} 
\newcommand{\Ssm}{Sierpi\'nski set in the sense of measure} 
\newcommand{\Sssm}{Sierpi\'nski sets in the sense of measure} 

\newcommand{\iA}{\mathcal{A}} 
\newcommand{\iH}{\mathcal{H}}

\newcommand{\iN}{\mathcal{N}}
\newcommand{\Rn}{$\RR^n$}

\title{Measurable Envelopes, Hausdorff Measures and Sierpi\'nski Sets}

\author{M\'arton Elekes \thanks{Partially supported by Hungarian Scientific
Foundation grant no. 37758 and F 43620}\\ 
\\Department of Analysis, E\"otv\"os Lor\'and University\\
Budapest, P\'azm\'any P\'eter s\'et\'any 1/c, 1117, Hungary\\
e-mail: {\tt emarci@cs.elte.hu}}

\begin{document}

\maketitle 

\begin{abstract}
We show that the existence of measurable envelopes of all subsets of $\RR^n$
with respect to the $d$ dimensional Hausdorff measure $(0<d<n)$ is
independent of $ZFC$. We also investigate the consistency of the existence
of $\Hd$-measurable \Sss. 
\end{abstract}

\insert\footins{\footnotesize{MSC codes: Primary 28A78; Secondary 03E35,
28E15}}
\insert\footins{\footnotesize{Key Words: measurable envelope, Hausdorff
measure, Sierpi\'nski set, independent}}

\section*{Introduction}
The following definition was motivated by the theory of analytic sets.

\bd
Let $\iA$ be a $\sigma$-algebra of subsets of a set $X$. We call a set
$H\subset X$ \textit{small} (with respect to $\iA$) if
every subset of $H$ belongs to $\iA$. The $\sigma$-ideal of small subsets
is denoted by $\iA_0$. We say that $A\in \iA$ is a \textsl{measurable
envelope} of $H\subset X$ (with respect to $\iA$) if $H\subset A$ and for
every $B\in \iA$ such that $H\subset B\subset A$ we have $A\setminus B\in
\iA_0$.

I have learnt the terminology `every subset of $X$ has a
measurable envelope' from D. Fremlin. Another usual one is `$(X,\iA)$ admits
covers' (see e.g. \cite{Ke}), and `measurable hull w.r.t. $\iA$' is also used. 
\ed

For example it is not hard to see that if $\iA$ is the Borel, Lebesgue or
Baire $\sigma$-algebra in $\RR^n$, then $\iA_0$ is the $\sigma$-ideal of
countable, Lebesgue negligible and first category sets, respectively. One can
also prove that with respect to the Lebesgue or Baire $\sigma$-algebra, every
subset of $\RR^n$ has a measurable envelope, while in the case of the Borel
sets this is not true.
What makes these notions interesting is a theorem of Szpilrajn-Marczewski,
asserting that if every subset of $X$ has a measurable envelope, then $\iA$
is closed under the Souslin operation (see \cite[29.13]{Ke}). 

This problem has been considered for various $\sigma$-algebras for a long
time (see e.g. \cite{Ma} and \cite{Pa}). 

In our paper we investigate the case $\iA=\iA_\mu$, where $\mu$ is an outer
measure on $X$ and $\iA_\mu$ is the $\sigma$-algebra of $\mu$-measurable
sets (in the sense of Carath\'eodory). To the best of our knowledge this 
problem was posed by M. Laczkovich.

It is well known and trivial that in the $\sigma$-finite
case every subset has a measurable envelope. Therefore we turn to the 
Hausdorff measures, which are
probably the most natural examples of non-$\sigma$-finite measures. We prove,
that this question cannot be answered in $ZFC$.  

As an application we give a short proof of the known statement that the
existence of an $\iH^1$-measurable \Ss \ (see the definition below) is
consistent with $ZFC$ (this is proven for the so called `one dimensional
measures' in \cite[3.11]{DP}).

Finally, we investigate the existence of two kinds of \Sss\ measurable with
respect to Hausdorff measures:

\bd
A set $S\subset\RR^2$ is a \textsl{\Ssm} if $S$ is (one dimensional)
Lebesgue negligible on each vertical line, but co-negligible (that is, the
complement  of $S$ is negligible) on each horizontal line.
A set $S\subset\RR^2$ is a \textsl{\Ssc} if $S$ is countable on each
vertical line, but co-countable on each horizontal line.
\ed

On one hand, we prove that for $0<d<2$ the existence of
$\iH^d$-measurable \Sssm\ is independent of $ZFC$. (In the remaining cases the
answer is trivial.) On the other hand, we show in $ZFC$ that in the
non-trivial cases ($0<d\le 2$) there exists no $\iH^d$-measurable \Ssc.

\smallskip

\br
Instead of considering Hausdorff measures, it would also be natural to look at
non-$\sigma$-finite outer
measures on arbitrary sets in general. As the following example shows, there
is a $ZFC$ example here.

\bex
(\cite{Fr}) Put $X=\om_2\times\om_2$ and let $\nu$ be the outer measure on $X$ that is
0 for countable subsets and 1 otherwise. Define
$$
\mu(H)=\sum_{\al\in\om_2}
\left[\nu(H\cap(\{\al\}\times\om_2))+\nu(H\cap(\om_2\times\{\al\}))\right]; 
$$
that is, let $\mu(H)$ be the number of uncountable horizontal and vertical
sections. Than one can check that $\om_1\times\om_2$ has no measurable
envelope.
\eex


\er

\section{Measurable envelopes with respect to Ha\-us\-dorff measures}

Let $\Hd$ denote the $d$-dimensional Hausdorff measure in $\RR^n$. Instead of
`measurable envelope w.r.t.
$\iA_{\Hd}$' we will write `measurable envelope w.r.t.
$\Hd$'.

If $d=0$ or $d>n$ then every subset of $\RR^n$ is $\Hd$-measurable, hence
every subset has a measurable envelope. If $d=n$ then $\Hd$ is $\sigma$-finite,
therefore we get the same conclusion. In the remaining cases we have

\bt\lab{env}
The following statement is independent of $ZFC$: for all $n\in\NN$ and $0<d<n$
every subset of $\RR^n$ has a measurable envelope with respect to $\Hd$.
\et

Before the proof we need two lemmas. ($\lambda$ denotes the one-dimensional
Lebesgue measure here.)

\bl\lab{isom}
Let $B\subset\RR^n$ be Borel such that $0<\Hd(B)<\infty$. Then there exists
a bijection $f$ between $B$ and the interval $I=[0,\Hd(B)]$ such that both $f$
and its inverse preserve Borel sets and measurable sets, and
which is measure preserving between the measure spaces $(B,\Hd)$ and
$(I,\lambda)$.
\el

\bp
\cite[12.B]{Ke} and \cite[17.41]{Ke}.
\ep

Now we turn to the second lemma.

\bd\lab{adddef}
$\textsl{add}\ \iN$ is the minimal cardinal $\ka$ for which there are $\ka$
Lebesgue negligible sets $A_\al\ (\al<\ka)$ such that $\cup_{\al<\ka}A_\al$
is of positive outer measure.
\ed

\br\lab{unio}
Note that if the sets $A_\al\ (\al<\la)$ are Lebesgue measurable for some
$\la<\textsl{add}\ \iN$, then so is their union $\cup_{\al<\la}A_\al$, as it
can be shown
by well-ordering the sets and noting that all but countably many of them
must be almost covered by the preceding sets.
\er

The following lemma is essentially contained in \cite[2.5.10]{Fe}.

\bl\lab{loc}
Let $0<d<n$ and suppose $\textsl{add}\ \iN=2^\om$. Then there exists a
disjoint family $\{M_\al:\al<2^\om\}$ of $\Hd$-measurable subsets of
$\RR^n$ of finite $\Hd$-measure, such that a set $H\subset\RR^n$ is
$\Hd$-measurable iff $H\cap M_\al$  is $\Hd$-measurable for every
$\al<2^\om$.
\el

\bp
Let $\{B_\al:\al<2^\om\}$ be be an enumeration of the Borel subsets of \Rn\
 of finite $\Hd$-measure, and put
$M_\al=B_\al\setminus(\cup_{\beta<\al}B_\beta)$. These are clearly pairwise
disjoint sets of finite $\Hd$-measure. Moreover $\textsl{add}\ \iN=2^\om$
together with the above remark and Lemma \ref{isom} applied to $B_\al $
gives that $M_\al$ is
$\Hd$-measurable for every $\al<2^\om$. The other direction being trivial
we only have to verify that if $H\subset\RR^n$ is such that $H\cap M_\al$
is $\Hd$-measurable for every $\al<2^\om$, then $H$ is itself
$\Hd$-measurable. Let $A\subset\RR^n$ arbitrary. We show that
$$
\Hd(A)\geq\Hd(H\cap A)+\Hd(H^C\cap A).
$$
We can obviously assume that $\Hd(A)<\infty$ and thus we can find a Borel set $B$ such that
$A\subset B$ and $\Hd(A)=\Hd(B)$, therefore $B=B_\al$ for some
$\al<2^\om$. Thus it is sufficient to prove that $H\cap B_\al$ is
$\Hd$-measurable, as that would imply 
$$
\Hd(A)=\Hd((H\cap B_\al)\cap A)+\Hd((H\cap B_\al)^C\cap A),
$$
but $(H\cap B_\al)\cap A=H\cap A$ and $(H\cap B_\al)^C\cap A=H^C\cap A$. In order to show that $H\cap B_\al$ is $\Hd$-measurable,
note that $B_\al=\cup_{\beta\le\al}M_\beta$, therefore $H\cap
B_\al=\cup_{\beta\le\al}(H\cap M_\al)$, which is again easily seen to be
$\Hd$-measurable.
\ep

\bd
$\textsl{non}^*\iN$ is the minimal cardinal $\ka$ such that in every
subset of the reals of positive outer Lebesgue measure we can find a subset
of positive outer Lebesgue measure and of cardinal $\le\ka$.

$\textsl{cov}\ \iN$ is the minimal cardinal $\ka$ such that $\RR$ can be
covered by $\ka$ Lebesgue negligible sets.
\ed

\br
$\textsl{non}^*\iN<\textsl{cov}\ \iN$ is consistent with $ZFC$ as it holds
in the so called `random real model', see \cite[Lemma 8]{LM}.
\er

Now we can turn to the proof of Theorem \ref{env}.

\bp
First we show that $\textsl{add}\ \iN=2^\om$ implies that for every $n\in\NN$
and $0<d<n$ every subset of $\RR^n$ has a measurable envelope with respect
to $\Hd$. (This proves that this statement is consistent, as
$\textsl{add}\ \iN=2^\om$ follows e.g. from $CH$ or $MA$.) Fix $n,\ d$ and
$H\subset\RR^n$. Let $\{M_\al:\al<2^\om\}$ be as in Lemma \ref{loc}. As
$M_\al$ is of finite measure for every $\al<2^\om$, we can find a
$\Hd$-measurable set $H_\al$ such that $H\cap M_\al\subset H_\al\subset
M_\al$ and $\Hd(H\cap M_\al)=\Hd(H_\al)$. We claim that
$A=\cup_{\al<2^\om}H_\al$ is a measurable envelope of $H$. Clearly $A$ is
$\Hd$-measurable by Lemma \ref{loc}. Suppose $H\subset B\subset A$, $B$ is
$\Hd$-measurable and $C\subset A\setminus B$. We want to show that $C$ is
measurable, therefore it is sufficient to check that $C\cap M_\al$ is
$\Hd$-measurable for every $\al<2^\om$, which is obvious, as it is of
$\Hd$-measure zero.

Next we prove that for $n=2$ and $d=1$ it is consistent that there
exists a subset of the plane without a $\iH^1$-measurable envelope. We
assume $\textsl{non}^*\iN<\textsl{cov}\ \iN$. One can easily find a set
$A\subset \RR$ of full outer measure and of cardinal $\textsl{non}^*\iN$,
and we claim that $A\times\RR$ has no $\iH^1$-measurable
envelope. Otherwise, if $M$ is such an envelope, then it is (one
dimensional) Lebesgue measurable on each vertical and horizontal line,
therefore it is Lebesgue negligible on all vertical lines over
$\RR\setminus A$ and co-negligible on all horizontal lines. As
$\textsl{non}^*\iN<\textsl{cov}\ \iN$, $\RR\setminus A$ is not negligible,
hence we can choose a set $B\subset\RR\setminus A$ of positive outer
measure and of cardinal $\textsl{non}^*\iN$. Then the projection of the set
$(B\times\RR)\cap M$ to the second coordinate consists of
$\textsl{non}^*\iN$ zero sets, on the other hand, it is the whole line, a
contradiction. 
\ep

\br The second direction of this proof (the last paragraph, in which we show a
set without a measurable envelope)
is due to D. Fremlin (\cite{Fr}).
In fact, it is not much harder to see that
$\textsl{non}^*\iN<\textsl{cov}\ \iN$ implies the existence of subsets of
\Rn\  without $\Hd$-measurable envelopes for any $0<d\le
[\frac{n}{2}]$ and $n\ge2$. Indeed, we can replace
$\RR\times\RR$ by the square of a $d$ dimensional Cantor set in
$\RR^{[\frac{n}{2}]}$ of positive and finite $\Hd$-measure, and
repeat the above argument.
\er

However, we do not know the answer to the following question.

\bq
Is it consistent that there exists a subset of \Rn\  without
a $\Hd$-measurable envelope for $n=1$, $0<d<1$ and for $n\ge 2$,
$[\frac{n}{2}]<d<n$?
\eq

\section{Hausdorff measurable \Sss}

As all subsets of the plane are $\iH^d$-measurable for $d=0$ or $d>2$, the
existence of $\iH^d$-measurable \Sssc\  for these $d$-s is equivalent to the
existence  of \Sssc, which is known to be equivalent to $CH$ (see \cite{Tr}). The following theorem answers the question for the other $d$-s.

\bt\lab{card}
For $0<d\le 2$ there exists no $\iH^d$-measurable \Ssc.
\et

\bp
For $d=2$ the statement is obvious by the Fubini Theorem. Let $0<d<2$ and
$C_{\frac{d}{2}}$ be a symmetric self-similar Cantor set in $[0,1]$ of
dimension $\frac{d}{2}$.

\bl
There exists $0<c<\infty$, such that
$$
\Hd|C_{\frac{d}{2}}\times C_{\frac{d}{2}}=
c\ (\iH^{\frac{d}{2}}|C_{\frac{d}{2}} \times
\iH^{\frac{d}{2}}|C_{\frac{d}{2}}).
$$
\el

\bp
Since $K=C_{\frac{d}{2}}\times C_{\frac{d}{2}}$ is also self-similar, it is
easy to see that $0<\Hd(K)<\infty$, therefore if we let 
$$
c=\frac{\Hd(K)}{(\iH^{\frac{d}{2}}(C_{\frac{d}{2}}))^2},
$$
then $0<c<\infty$ and the two above outer measures agree on $K$. By the
self-similarity of $K$ they also agree on the basic open sets, and as any
open set is the disjoint union of countably many of these, the outer
measures agree on all open sets. Hence (by finiteness) on all Borel sets as
well, from which the lemma follows.
\ep

Now we can complete the proof of Theorem \ref{card} as follows. Note that
if $S$ is an $\iH^d$-measurable \Ssc, then $S\cap K$ is a
$\iH^d$-measurable set, which is countable on each vertical section of $K$
and co-countable on each horizontal section of $K$. But this gives a
contradiction, once we apply the previous lemma and the Fubini Theorem.
\ep

The question concerning \Sssm\  is more complicated. Just as above, for
$d=0$ or $d>2$ \  $\iH^d$-measurability is not really a restriction, and we
know that the existence of \Sssm\ is independent of $ZFC$ (see \cite[Theorem
2]{La}). For $d=2$ no such set can be $\iH^d$-measurable by the Fubini
Theorem, while in the remaining cases we have the following.

\bt
For $0<d<2$ the existence of $\iH^d$-measurable \Sssm\ is independent of
$ZFC$. 
\et

\bp
On one hand, for example $\textsl{non}^*\iN<\textsl{cov}\ \iN$ implies that
there are no \Sss\  of any kind (\cite[Theorem 2]{La}).

On the other hand, we assume $\textsl{add}\ \iN=2^\om$ and prove the
existence of $\iH^d$-measurable \Sssm, separately for all $0<d<1$, $d=1$
and all $1<d<2$.

If $d=1$, then our statement is a consequence of
\cite[3.10]{DP}, but we present another proof here. By \cite[Theorem
2]{La} and $\textsl{add}\ \iN=2^\om$ we can find a \Ssm, and by \ref{env}
this set has a $\iH^d$-measurable envelope. It is not hard to check, that
this envelope possesses the required properties. 

Now let $0<d<1$. Enumerate the Borel subsets of $\RR^2$ of positive finite
$\iH^d$-measure as $\{B_\al:\al<2^\om\}$ and also $\RR$ as
$\{x_\al:\al<2^\om\}$. We can assume, that $S$ is a \Ssm, such that the
cardinality of every vertical section is less than $2^\om$ (the proof in
\cite{La} provides such a set).
Then put
$$
S_1=S\cup\bigcup_{\al<2^\om}\left[B_\al\setminus\left(\bigcup_{\beta<\al}B_\beta\cup(\{x_\beta:\beta<\al\}\times\RR)\right)\right].
$$
$S_1$ is a \Ssm\ as its horizontal sections contain the horizontal sections
of $S$, the vertical section over $x$ is still of (one dimensional) Lebesgue
measure zero,
since $\textsl{add}\ \iN=2^\om$, and $x=x_\al$ for some $\al<2^\om$ so this
section is increased only in the first $\al$ steps, and always by a set of
finite $\iH^d$-measure, therefore of zero Lebesgue measure. What remains to
check is that $S_1$ is $\iH^d$-measurable. Like above, by the Borel
regularity of $\iH^d$ it is sufficient to show that
$$
\iH^d(B)=\iH^d(S_1\cap B)+\iH^d(S_1^C\cap B)
$$
holds for every Borel set $B$ of finite and positive $\iH^d$-measure,
therefore we only have to prove that $S_1\cap B_\al$ is $\iH^d$-measurable
for every $\al<2^\om$. We show this by induction on $\al$ as follows. Put
$$
A_\al=\bigcup_{\beta<\al}B_\beta\cup(\{x_\beta:\beta<\al\}\times\RR).
$$
Then
$$
S_1\cap B_\al=S_1\cap
[(B_\al\setminus A_\al) \cup (B_\al\cap A_\al)]=
$$
$$
[S_1\cap(B_\al\setminus A_\al)] \cup [S_1\cap(B_\al\cap A_\al)]=
[B_\al\setminus A_\al] \cup [S_1\cap(B_\al\cap A_\al)]=
$$
$$
[B_\al\setminus A_\al] \cup 
\left[\bigcup_{\beta<\al}(B_\al\cap B_\beta\cap S_1)\right] \cup 
\left[\bigcup_{\beta<\al}(B_\al\cap(\{x_\beta\}\times\RR)\cap S_1)\right].
$$
Now we apply Lemma \ref{isom} to $B_\al$ in view of $\textsl{add}\
\iN=2^\om$. Then the first expression in the
last line is clearly $\iH^d$-measurable (we may apply
Lemma \ref{isom} and the Remark following Definition \ref{adddef} to $B_\al$), and the
same conclusion holds for the second
expression by our inductional hypotheses. In order to check measurability for
the last one we note that $B_\al\cap(\{x_\beta\}\times\RR)\cap S$ is of
cardinal less than $2^\om$ for every $\beta<\al$, thus $\iH^d$-negligible,
but when we construct
$B_\al\cap(\{x_\beta\}\times\RR)\cap S_1$ out of this set, we increase it
only in the first $\beta$ steps, and always by a $\iH^d$-measurable set.

Finally, let $1<d<2$. As above, let $\{B_\al:\al<2^\om\}=\{B\subset\RR: B
\textrm{ Borel},\iH^d(B)<\infty\}$ and also $\{x_\al:\al<2^\om\}=\RR$. Put
$$
D_\al=B_\al\setminus[\bigcup_{\beta<\al}B_\beta \cup
(\{x_\beta:\beta<\al\}\times\RR) \cup (\RR\times\{x_\beta:\beta<\al\})]
$$
for every $\al<2^\om$.
Since $D_\al\subset B_\al$,\ \  $D_\al$ is (two dimensional) Lebesgue
negligible, therefore $\{x\in\RR: \la((\{x\}\times\RR)\cap D_\al)>0\}$ is
Lebesgue negligible, thus contained in a Borel set $N_\al$ of Lebesgue
measure zero. Define
$$
S_1= \left(  S\cup\left[\bigcup_{\al<2^\om}(D_\al\setminus(N_\al\times\RR))\right]  \right)  \setminus
\left[\bigcup_{\al<2^\om}(D_\al\cap(N_\al\times\RR))\right].
$$
First we check that $S_1$ is a \Ssm. If $\RR\times\{x\}$ is a horizontal
line, then $x=x_\al$ for some $\al<2^\om$. $D_\xi$ and $\RR\times\{x_\al\}$
are disjoint for every $\xi>\al$, therefore our set is not modified after
the first $\al$ steps. Hence
$$
[S\setminus S_1]\cap(\RR\times\{x_\al\}) \subset
\left[\bigcup_{\xi\le\al}(D_\xi\cap(N_\xi\times\RR))\right]\cap(\RR\times\{x_\al\}),
$$
which is Lebesgue negligible on this horizontal line by $\textsl{add}\
\iN=2^\om$. Similarly, on a vertical line $\{x_\al\}\times\RR$
$$
S_1\cap(\{x_\al\}\times\RR) \subset
\left[S\cup\bigcup_{\xi\le\al}(D_\xi\setminus(N_\xi\times\RR))\right] \cap
(\{x_\al\}\times\RR),
$$
which is again a zero set. What remains to show is the $\iH^d$-measurability
of $S_1$. It is again sufficient to prove by induction on $\al$ that
$S_1\cap B_\al$ is $\iH^d$-measurable for every $\al<2^\om$. Just like
above, we apply Lemma \ref{isom} to $B_\al$.
$$
S_1\cap B_\al=
$$
$$
[S_1\cap D_\al] \cup
$$
$$
\left[S_1\cap B_\al\cap
\left(\bigcup_{\beta<\al}B_\beta \cup
(\{x_\beta:\beta<\al\}\times\RR)\cup
(\RR\times\{x_\beta:\beta<\al\})\right)\right],
$$
where the expression in the first brackets equals
$D_\al\setminus(N_\al\times\RR)$, which is clearly $\iH^d$-measurable ($N_\al$ is Borel and to $D_\al$ we can apply the usual argument in
$B_\al$). As for the $\iH^d$-measurability of the second expression, it equals
$$
B_\al\cap
$$
$$
\left([\bigcup_{\beta<\al}(B_\beta\cap S_1)] \cup 
[(\{x_\beta:\beta<\al\}\times\RR)\cap S_1] \cup
[(\RR\times\{x_\beta:\beta<\al\})\cap S_1]\right).
$$
We apply the usual argument in $B_\al$.
$B_\beta\cap S_1$ is
$\iH^d$-measurable for every $\beta<\al$ by the inductional hypothesis.
Moreover, as $d>1$, $\{x_\beta:\beta<\al\}\times\RR$ as well as
$\RR\times\{x_\beta:\beta<\al\}$ is of $\iH^d$-measure zero in $B_\al$ for every
$\beta<\al$. Therefore the proof is complete.
\ep

However, we do not know the answer to the following.

\bq
Let $0<d<2$. Is it consistent, that there exists a \Ssm, but it cannot be
$\iH^d$-measurable?
\eq

\ba
I am greatly indebted to my advisor Professor Mikl\'os Lacz\-ko\-vich for the many
useful discussions.

\end{document}